\documentclass[12pt,a4paper]{amsart}
\usepackage[utf8]{inputenc}	
\usepackage[T1]{fontenc}
\usepackage[in,headings]{fullpage}
\usepackage{amsbsy, amscd, amsfonts, amsmath, amsrefs, amssymb, amsthm}
\usepackage{graphicx}
\usepackage{float}
\usepackage{color}   
\usepackage[colorlinks=true,linkcolor=blue,citecolor=blue,urlcolor=blue]{hyperref}
\usepackage{comment}
\includecomment{A}

\newtheorem{theorem}{Theorem}

\newtheorem{proposition}[theorem]{Proposition}
\newtheorem{lemma}[theorem]{Lemma}

\theoremstyle{definition}

\theoremstyle{remark}

\newcommand{\R}{\mathbf{R}}

\renewcommand{\Re}{\mathop{\mathrm{Re}}\nolimits}

\newcommand{\Rzeta}{\mathop{\mathcal R }\nolimits}

\newfont{\cmbsy}{cmbsy10}
\newfont{\cmmib}{cmmib10}

\DeclareMathOperator{\Res}{Res}

\begin{document}

\title[Integral representation due to Gabcke]{An Integral representation of $\Rzeta(s)$ due to Gabcke.}
\author[Arias de Reyna]{J. Arias de Reyna}
\address{%
Universidad de Sevilla \\ 
Facultad de Matem\'aticas \\ 
c/Tarfia, sn \\ 
41012-Sevilla \\ 
Spain.} 




\email{arias@us.es, ariasdereyna1947@gmail.com}


\begin{abstract}
Gabcke proved a new integral expression for the auxiliary Riemann function
\[\Rzeta(s)=2^{s/2}\pi^{s/2}e^{\pi i(s-1)/4}\int_{-\frac12\searrow\frac12}
\frac{e^{-\pi i u^2/2+\pi i u}}{2i\cos\pi u}U(s-\tfrac12,\sqrt{2\pi}e^{\pi i/4}u)\,du,\]
where $U(\nu,z)$ is the usual parabolic cylinder function. 

We give a new, shorter proof, which avoids the use of the Mordell integral.  And we write it in the form
\begin{equation}
\Rzeta(s)=-2^s \pi^{s/2}e^{\pi i s/4}\int_{-\infty}^\infty \frac{e^{-\pi x^2}H_{-s}(x\sqrt{\pi})}{1+e^{-2\pi\omega x}}\,dx.
\end{equation}
where $H_\nu(z)$ is the generalized Hermite polynomial.
\end{abstract}

\maketitle
\section{Introduction}
The auxiliary function of Riemann is defined by the integral
\[\Rzeta(s)=\int_{0\swarrow1}\frac{x^{-s} e^{\pi i x^2}}{e^{\pi i x}-e^{-\pi i x}}\,dx.\]
The position of the zeros of this  function is connected with the zeros of the Riemann zeta function \cite{A166}. In a paper by Gabcke \cite{G}, posted on arXiv and not published, he proved an integral 
expression for $\Rzeta(s)$ using the parabolic cylinder function $U(\nu,z)$ as defined in \cite{AS}.
We prefer to give it in terms of the Hermite functions \cite{L}. We show that 
\begin{equation}\label{main}
\Rzeta(s)=-2^s \pi^{s/2}e^{\pi i s/4}\int_{-\infty}^\infty \frac{e^{-\pi x^2}H_{-s}(x\sqrt{\pi})}{1+e^{-2\pi\omega x}}\,dx,\qquad \omega=e^{\pi i/4}.
\end{equation}
We give a different proof of that of Gabcke, which is simpler and does not use the Mordell integral. Furthermore, our proof presents some integrals that may be of independent interest.

\section{Somme Lemmas}

\begin{lemma}\label{L:sinc}
For $z=x+iy$ with $0<x-y<1$, we have
\[\int_L\frac{e^{-z\zeta}}{1+e^{-\zeta}}\,d\zeta=\frac{\pi}{\sin\pi z},\]
where $L$ is a line through $0$ in direction $e^{\pi i/ 4}$.
\end{lemma}
\begin{proof} This is an exercise in complex analysis. 
Let  $\omega=e^{\pi i/4}$, then $\zeta=\omega u$ and 
\[\int_L\frac{e^{-z\zeta}}{1+e^{-\zeta}}\,d\zeta=
\omega\int_{-\infty}^\infty \frac{e^{-\omega z u}}{1+e^{-\omega u}}\,du.\]
Since  $\omega z=2^{-1/2}(1+i)(x+iy)$ our hypothesis implies that $0<\alpha:=\Re(\omega z)<2^{-1/2}$ the integrand decays as $e^{-(\alpha-2^{-1/2})u}$ for $u\to-\infty$, and as $e^{-\alpha u}$ when $u\to+\infty$. Therefore, the integral is well defined and absolutely convergent.
It is clear that the same happens to the integral along $L+w$ for any complex number $w$. We have 
\[\int_L\frac{e^{-z\zeta}}{1+e^{-\zeta}}\,d\zeta-\int_{L+2\pi i}\frac{e^{-z\zeta}}{1+e^{-\zeta}}\,d\zeta=2\pi i \Res_{\zeta=\pi i}\frac{e^{-z\zeta}}{1+e^{-\zeta}}=2\pi i e^{-\pi i z}.\]
On the other hand, by the periodicity of $e^z$ we have
\[\int_{L+2\pi i}\frac{e^{-z\zeta}}{1+e^{-\zeta}}\,d\zeta=
\int_{L}\frac{e^{-z(\zeta+2\pi i)}}{1+e^{-\zeta}}\,d\zeta=e^{-2\pi i z}\int_L\frac{e^{-z\zeta}}{1+e^{-\zeta}}\,d\zeta.\]
It follows that 
\[\int_L\frac{e^{-z\zeta}}{1+e^{-\zeta}}\,d\zeta=2\pi i \frac{e^{-\pi i z}}{1-e^{-2\pi i z}}=\frac{\pi}{\sin\pi z}.\qedhere\]
\end{proof}

The parabolic cylinder functions are solutions of the differential equation
\[\frac{d^2w}{dz^2}=\Bigl(\frac{z^2}{4}+a\Bigr) w.\]
The principal solution $U(a,z)$ is defined by the condition $U(a,z)\sim z^{-a-\frac12}e^{-z^2/4}$ it is an entire function in $z$ and in $a$. There are two other alternative notation, Whittaker's $D_\nu(z)$ and Hermite functions $H_\nu(z)$. They are related by 
\[D_\nu(z)=U(-\nu-1/2,z)=2^{-\nu/2}e^{-z^2/4}H_\nu(z/\sqrt{2}).\]
Therefore, they are equivalent. We will use the notation of Hermite functions that is more familiar to a general mathematician. It is true that the Hermite functions are the most unusual in the specific bibliography. One exception is the book by Lebedev \cite{L}. For $\nu$ a nonnegative integer the functions $H_\nu(z)$ are the usual Hermite polynomials. Since they are entire functions, we may define them by its power series expansion \cite{L}*{eq.~(10.4.3)} (for $\nu$ not a nonnegative integer)
\[H_\nu(z)=\frac{1}{2\Gamma(-\nu)}\sum_{n=0}^\infty \Gamma\Bigl(\frac{n-\nu}{2}\Bigr)\frac{(-2z)^n}{n!}.\]
In Lebedev \cite{L}*{eq.~(10.5.2)} we find the following integral representation
\begin{equation}\label{L:Int}
H_\nu(z)=\frac{1}{\Gamma(-\nu)}\int_0^\infty e^{-t^2-2tz}t^{-\nu-1}\,dt,\qquad \Re\nu<0.\end{equation}

\begin{proposition}\label{P:IntH}
For any complex numbers $s$ and $z$ we have 
\begin{equation}
\int_{0\uparrow} x^s e^{x^2/4-x z}\,dx=2i\sqrt{\pi} e^{-z^2}H_s(z),
\end{equation}
where we integrate along any line of direction $e^{i\theta}$ with $\pi/4<\theta<3\pi/4$ that cuts the positive real axis. 
\end{proposition}
\begin{proof}
An application of Cauchy's Theorem proves that the integral do not depend on the line of integration in the conditions of the proposition. We take a vertical line in the proof. First, we prove that the integral is an entire function of $z$ and of $s$. 
Consider the line of integration $x=1+i u$ with $u\in\R$. Assume that $|s|\le R$ and $|z|\le R$ for some positive real number $R>1$. We have 
\[|x^{s}|\le|\exp((\sigma+it)(\log|1+i u|+i\arctan u))|\le |1+i u|^{\sigma} e^{\pi|t|}\le e^{\pi R}(1+|u|)^R.\]
\[x^2/4- zx=\tfrac14 (1+i u)^2- z-izu=\tfrac14-\tfrac14u^2+\tfrac12iu-z-izu,\]
so that
\[|\exp(x^2/4- zx)|\le e^{-u^2/4+R+R|u|+\frac14}.\]
Therefore,
\[|x^{-s} e^{x^2/4- zx}|\ll_R(1+|u|)^R e^{-u^2/4+R|u|}.\]
Hence, the integrand is bounded by an absolutely integrable function  of $u$.  Since the integrand is holomorphic in $s$ and $z$, it follows that the integral defines an entire function of $z$ and $s$. 

When $\sigma>-1$, applying Cauchy's Theorem, we may move the line of integration to $x=iu$, since  the singularity at $x=0$ is integrable, and the bound we have proved slightly modified still applies. 
We find that our integral is equal in this case to 
\[i e^{\pi i s/2}\int_0^\infty  u^{s}e^{-u^2/4-iz u}\,du+
ie^{-\pi i s/2}\int_0^\infty u^{s}e^{-u^2/4+izu}\,du.\]
By \eqref{L:Int} this is equal to
\[ie^{\pi i s/2}2^{1+s}\Gamma(1+s)H_{-s-1}(iz)+ie^{-\pi i s/2}2^{1+s}\Gamma(1+s)H_{-s-1}(-iz).\]
Equivalently
\[i2^{s+1}\Gamma(1+s)\bigl[e^{\pi i s/2}H_{-s-1}(iz)+e^{-\pi i s/2}H_{-s-1}(-iz)\bigr].\]
In Lebedev \cite{L}*{eq.~(10.3.4)} it is proved that
\[H_\nu(z)=\frac{2^\nu\Gamma(\nu+1)}{\sqrt{\pi}}e^{z^2}\bigl[e^{\pi i \nu/2}H_{-\nu-1}(iz)+e^{-\pi i \nu/2}H_{-\nu-1}(-iz)\bigr].\]
Hence, our integral 
\[\int_{0\uparrow} x^s e^{x^2/4-x z}\,dx=
2i\sqrt{\pi}e^{-z^2}H_s(z),\]
as we wanted to prove. This has been proved only for $\sigma>-1$ but is then true for any $s$ and $z$ by analyticity.
\end{proof}

\section{Main Theorem}

\begin{theorem}
For any complex number $s$ we have
\begin{equation}
\Rzeta(s)=-2^s \pi^{s/2}e^{\pi i s/4}\int_{-\infty}^\infty \frac{e^{-\pi x^2}H_{-s}(x\sqrt{\pi})}{1+e^{-2\pi\omega x}}\,dx.
\end{equation}
\end{theorem}
\begin{proof}
By definition
\[\Rzeta(s)=\int_{0\swarrow1}\frac{x^{-s}e^{\pi i x^2}}{e^{\pi i x}-e^{-\pi i x}}\,dx.\]
Take $x=\frac12-\omega u$ as the line of integration. Any $x$ in this line satisfies the conditions of lemma \ref{L:sinc}. Therefore,
\[\Rzeta(s)=\frac{1}{2\pi i}\int_{0\swarrow1}x^{-s}e^{\pi i x^2}\Bigl(\int_L\frac{e^{-x\zeta}}{1+e^{-\zeta}}\,d\zeta\Bigr)\,dx,\]
where $L$ is the line $\zeta=\omega v$ with $v\in \R$.
The function 
\[-i\frac{x^{-s}e^{\pi i x^2-x\zeta}}{1+e^{-\zeta}},\qquad \zeta=\omega v,\quad x=\frac12-\omega u,\]
is absolutely integrable in $(u,v)\in\R^2$. In fact, we have for $|s|\le R$
\begin{align*}
|x^{-s}|&=|\exp(-(\sigma+it)(\log|\tfrac12-\omega u|+i\arg(\tfrac12-\omega u))|\\
&\ll_R\max(2^{\frac{3R}{2}},(1+|u|)^R)e^{\pi R}\ll_R(1+|u|)^R.
\end{align*}
\[\pi i x^2-x\zeta=\pi i(\tfrac12-\omega u)^2-(\tfrac12-\omega u)\omega v=
\tfrac14\pi i-\pi u^2-\pi i \omega u-\tfrac12\omega v+iuv,\]
so that
\[|\exp(\pi i x^2-x\zeta)|=\exp(-\pi u^2+\pi u/\sqrt{2}-v/2\sqrt{2}).\]
$1+e^{-\zeta}=1+e^{-\omega v}$ is a continuous function that does not vanish for $v\in\R$ and 
\[|1+e^{-\zeta}|=|1+e^{-v/\sqrt{2}-iv/\sqrt{2}}|\gg \max(e^{-v/\sqrt{2}}-1, 1).\]
Hence,
\[\Bigl|-i\frac{x^{-s}e^{\pi i x^2-x\zeta}}{1+e^{-\zeta}}\Bigr|\ll_R (1+|u|)^R\exp(-\pi u^2+\pi u/\sqrt{2})\cdot\frac{e^{-v/2\sqrt{2}}}{\max(e^{-v/\sqrt{2}}-1, 1)},\]
is bounded by the product of an integrable function of $u$ by an integrable function of $v$.

Hence, we may interchange the order of integration so that 
\[\Rzeta(s)=\frac{1}{2\pi i}\int_L\Bigl(\int_{0\swarrow1}x^{-s}e^{\pi i x^2-x\zeta}\,dx\Bigr)\frac{d\zeta}{1+e^{-\zeta}}.\]
In the inner integral change variable putting $\xi=2\sqrt{\pi}\omega x$, when $x$ runs through the line $0\swarrow1$ the variable $\xi$ runs through $\sqrt{\pi}\omega-2i\sqrt{\pi}u$ that is a vertical line $0\downarrow$ so that by Proposition \ref{P:IntH}
\begin{align*}
\int_{0\swarrow1}x^{-s}e^{\pi i x^2-x\zeta}\,dx=-\frac{2^{s}\pi^{s/2}e^{\pi i s/4}}{2\sqrt{\pi}\omega}\int_{0\uparrow}\xi^{-s}e^{\xi^2/4-\xi\zeta/2\sqrt{\pi}\omega}\,d\xi\\
=-\omega 2^{s}\pi^{s/2}e^{\pi i s/4}e^{-\frac{\zeta^2}{4\pi i}}H_{-s}\Bigl(\frac{\zeta}{2\sqrt{\pi}\omega}\Bigr).
\end{align*}
Hence, 
\[\Rzeta(s)=-\frac{\omega}{2\pi i} 2^{s}\pi^{s/2}e^{\pi i s/4}\int_{L}\frac{e^{-\frac{\zeta^2}{4\pi i}}H_{-s}(\frac{\zeta}{2\sqrt{\pi}\omega})}{1+e^{-\zeta}}\,d\zeta.\]
Finally, take the parametrization $\zeta=2\pi \omega v$ of the integration line to get
\[\Rzeta(s)=-2^{s}\pi^{s/2}e^{\pi i s/4}\int_{-\infty}^\infty\frac{e^{-\pi v^2}H_{-s}(\sqrt{\pi} v)}{1+e^{-2\pi\omega v}}\,dv.\qedhere\]
\end{proof}

\begin{figure}[H]
\begin{center}
\includegraphics[width=0.6\hsize]{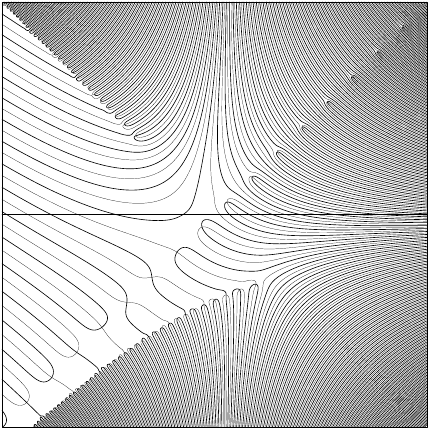}
\caption{x-ray of $H_{-\frac12-10i}(z\sqrt{\pi})e^{-\pi z^2}(1+e^{-2\pi\omega z})^{-1}$ in the square $(-6,6)^2$}
\end{center}
\end{figure}

\vspace{-1cm}

\begin{figure}[H]
\begin{center}
\includegraphics[width=0.6\hsize]{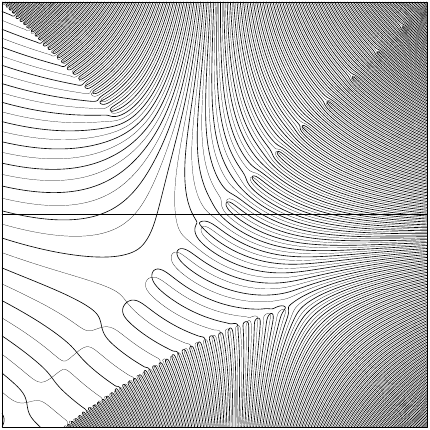}
\caption{x-ray of $H_{-\frac12-20i}(z\sqrt{\pi})e^{-\pi z^2}(1+e^{-2\pi\omega z})^{-1}$ in the square $(-6,6)^2$}
\end{center}
\end{figure}

\end{document}